\documentclass [12pt] {amsart}

\usepackage [T1] {fontenc}
\usepackage {amssymb}
\usepackage {amsmath}
\usepackage {amsthm}
\usepackage{systeme}
\usepackage{tikz}

\usepackage{url}

\usepackage{soul}

 \usepackage[titletoc,toc,title]{appendix}

\usetikzlibrary{arrows,decorations.pathmorphing,backgrounds,positioning,fit,petri}

\usetikzlibrary{matrix}

\usepackage{enumitem}
\usepackage{hyperref}

\DeclareMathOperator {\td} {td}

\DeclareMathOperator {\Th} {Th}

\DeclareMathOperator {\tp} {tp}

\DeclareMathOperator {\cl} {cl}

\DeclareMathOperator {\pr} {pr}

\DeclareMathOperator {\SL} {SL}

\DeclareMathOperator {\PGL} {PGL}
\DeclareMathOperator {\GL} {GL}

\DeclareMathOperator {\B} {B}

\DeclareMathOperator {\an} {an}

\DeclareMathOperator {\C} {\mathbb{C}}
\DeclareMathOperator {\R} {\mathbb{R}}
\DeclareMathOperator {\h} {\mathbb{H}}
\DeclareMathOperator {\Z} {\mathbb{Z}}
\DeclareMathOperator {\Q} {\mathbb{Q}}
\DeclareMathOperator {\A} {A}
\DeclareMathOperator {\ppr} {Pr}
\DeclareMathOperator {\reg} {reg}
\DeclareMathOperator {\F} {\mathbb{F}}
\DeclareMathOperator {\rar} {\rightarrow}
\DeclareMathOperator {\seq} {\subseteq}

\DeclareMathAlphabet\urwscr{U}{urwchancal}{m}{n}%
\DeclareMathAlphabet\rsfscr{U}{rsfso}{m}{n}
\DeclareMathAlphabet\euscr{U}{eus}{m}{n}
\DeclareFontEncoding{LS2}{}{}
\DeclareFontSubstitution{LS2}{stix}{m}{n}
\DeclareMathAlphabet\stixcal{LS2}{stixcal}{m} {n}

\theoremstyle {definition}
\newtheorem {definition}{Definition} [section]
\newtheorem{example} [definition] {Example}
\newtheorem* {claim} {Claim}
\newtheorem* {notation} {Notation}

\theoremstyle {plain}
\newtheorem {question} [definition] {Question}
\newtheorem {lemma} [definition] {Lemma}
\newtheorem {theorem} [definition] {Theorem}
\newtheorem {fact} [definition] {Fact}

\newtheorem {conjecture} [definition] {Conjecture}

\theoremstyle {remark}
\newtheorem {remark} [definition] {Remark}

\calclayout

\begin {document}

\title{Blurrings of the $j$-Function}

\author{Vahagn Aslanyan}
\address{Vahagn Aslanyan, School of Mathematics, University of East Anglia, Norwich, NR4 7TJ, UK}
\email{V.Aslanyan@uea.ac.uk}

\author{Jonathan Kirby}
\address{Jonathan Kirby, School of Mathematics, University of East Anglia, Norwich, NR4 7TJ, UK}
\email{Jonathan.Kirby@uea.ac.uk}

\thanks{Supported by EPSRC grant EP/S017313/1.}

\date{\today}

\keywords {Ax-Schanuel, Existential Closedness, $j$-function}

\subjclass[2010] {11F03, 03C60, 03C64}

\vspace*{-1cm}

\maketitle

\begin{abstract}
Inspired by the idea of blurring the exponential function, we define blurred variants of the $j$-function and its derivatives, where blurring is given by the action of a subgroup of $\GL_2(\C)$. For a dense subgroup (in the complex topology) we prove an Existential Closedness theorem which states that all systems of equations in terms of the corresponding blurred $j$ with derivatives have complex solutions,  except where there is a functional transcendence reason why they should not. For the $j$-function without derivatives we prove a stronger theorem, namely, Existential Closedness for $j$ blurred by the action of a subgroup which is dense in $\GL_2^+(\R)$, but not necessarily in $\GL_2(\C)$.

We also show that for a suitably chosen countable algebraically closed subfield $C \seq \C$, the complex field augmented with a predicate for the blurring of the $j$-function by $\GL_2(C)$ is model theoretically tame, in particular, $\omega$-stable and quasiminimal.
\end{abstract}

\section{Introduction}

Let $\h$ be the complex upper half-plane and let $j:\h \rightarrow \C$ be the modular $j$-function. It is invariant under the action of $\SL_2(\Z)$ on $\h$ and behaves nicely under the action of $\GL_2^+(\Q)$, namely, it satisfies certain algebraic ``functional equations'' given by the \emph{modular polynomials} (see Section \ref{subsect-modular-poly}). This gives rise to the Modular Schanuel conjecture (see \cite[Conjecture 8.3]{Pila-functional-transcendence}), henceforth referred to as MSC, which is a transcendence statement about the values of the $j$-function, and is an analogue of Schanuel's conjecture for the exponential function \cite[p. 30]{Lang-tr}. 

\begin{conjecture}[Modular Schanuel Conjecture]\label{MSC}
Let $z_1,\ldots,z_n \in \h$ be non-quadratic numbers with distinct $\GL_2^+(\Q)$-orbits. Then $$\td_{\Q}\Q(z_1,\ldots,z_n,j(z_1),\ldots,j(z_n))\geq n.$$
\end{conjecture}

This conjecture can be understood as a statement about certain ``overdetermined'' systems of polynomial equations of $2n$ variables not having solutions of the form $(z_1,\ldots,z_n, j(z_1),\ldots,j(z_n))$ with $z_k$'s having pairwise distinct $\GL_2^+(\Q)$-orbits. There is also a ``dual'' conjecture, known as \emph{Existential Closedness}, or concisely EC, stating that if such a system is not overdetermined, that is, having a solution does not contradict MSC, then there is a solution in $\C$. The notion of an overdetermined system of equations in $2n$ variables is captured by some properties of an algebraic variety $V\seq \C^{2n}$ known as $j$-\emph{broadness}, $j$-\emph{freeness}  {and $\h$-\emph{freeness}}. Roughly, $j$-broadness means that certain projections of $V$ are not too small, in particular, $\dim V \geq n$, while $j$-freeness means that no modular polynomial vanishes on the second $n$ coordinates of $V$,  {and $\h$-freeness means no $\GL_2^+(\Q)$-relation holds on the first $n$ coordinates of $V$}. We now formulate the EC conjecture postponing the precise definitions of $j$-broadness, $j$-freeness and $\h$-freeness to Section \ref{subsect-broad-free}.

{Throughout the paper algebraic varieties will be identified with the sets of their $\C$-points.} By a subvariety of $\h^n \times \C^n$ we mean a non-empty intersection of an algebraic variety in $\C^{2n}$ with $\h^n\times \C^n$. By abuse of notation we will let $j$ denote all Cartesian powers of itself. Similarly we let $\Gamma_j := \{ (\bar{z},j(\bar{z})): \bar{z}\in \h^n \}\seq \C^{2n}$ be the graph of $j$ in $\h^n \times \C^{n}$ for any $n$.

\begin{conjecture}[Existential Closedness for $j$]\label{j-ec}
Let $V \subseteq \h^n \times \C^n$ be an irreducible $j$-broad, $j$-free and  {$\h$-free} variety defined over $\C$. Then $V\cap \Gamma_j \neq \emptyset$.
\end{conjecture}

While MSC is considered out of reach, EC seems to be more tractable. In particular, EC was proven under the additional assumption that the projection of $V$ on the first $n$ coordinates has dimension $n$ (which in fact implies $j$-broadness and $\h$-freeness) by Eterovi\'c and Herrero in \cite{eterovic-herrero} (and some other related results were also obtained there), and a functional (differential) analogue of EC was established in \cite{Aslanyan-Eterovic-Kirby-Diff-EC-j} by Aslanyan, Eterovi\'c, Kirby. This paper is devoted to proving some weak versions of EC. 

There are many analogies between the $j$-function and the complex exponentiation, especially from a transcendence point of view. In particular, MSC is the modular analogue of the famous Schanuel conjecture for $\exp$, and the EC conjecture for $j$ is an analogue of a similar conjecture for $\exp$, known as \emph{Exponential Closedness}, posed by Zilber in the early 2000's \cite{Zilb-pseudoexp}. Moreover, the aforementioned theorems proven in \cite{eterovic-herrero} are analogous to some results on the exponential function established in \cite{brown-masser,daquino--fornasiero-terzo}, though the methods are quite different. Similarly, the work in this paper is related to previous work on the blurred complex exponentiation \cite{Kirby-blurred},  {but we also use some new tools here such as combining methods from complex and real analytic geometry and o-minimality.}

\begin{definition}
Given a subgroup $G \subseteq \GL_2(\C)$, let $\B_j^G\seq \C^2$ be the relation $\{ (z, j(gz)): g \in G, gz \in \h \}$. By abuse of notation we also let $\B^G_j$ denote the set $$\{ (z_1,\ldots,z_n, j(g_1z_1),\ldots, j(g_nz_n)):  g_k \in G, g_kz_k \in \h \mbox{ for all } k \}$$ for every $n$. We think of $\B_j^G$ as the graph of $j$ \emph{blurred} by the action of $G$.
\end{definition}

Let us consider a few examples.

\begin{example}
    \begin{itemize}
         \item[]
        \item When $G \seq \SL_2(\Z)$, we have $\B_j^G = \Gamma_j$.
        \item $\B_j^{\GL_2(\C)} = \C^2$.
        \item $\B_j^{\GL^+_2(\R)} = \h \times \C$.
    \end{itemize}
\end{example}

 {Observe that Conjecture \ref{j-ec} implies immediately that if $V$ is $j$-broad, $j$-free and $\h$-free then $V\cap \B_j^G \neq \emptyset$ for any group $G$.}

The following is one of our main results.

\begin{theorem}\label{j-blur-ec-intro}
If $V \seq \C^{2n}$ is a $j$-broad and $j$-free variety and $G$ is a dense subgroup of $\GL_2(\C)$ in the complex topology then $ V \cap \B_j^G$ is dense in $ V$, and hence it is non-empty.
\end{theorem}

Note that while in Conjecture \ref{j-ec} the variety is required to be $\h$-free, in the above theorem and in the other main results of the paper we do not have such a constraint. The difference is that in the conjecture we want the variety to intersect the graph of a function, and this can be used to show that $\h$-freeness is needed there. For example, the variety $V \seq \h^2 \times \C^2$ given by two equations $z_1=z_2, w_1=w_2+1$ is $j$-broad and $j$-free, but not $\h$-free, and it obviously does not intersect $\Gamma_j$.


We prove a more general theorem incorporating the (first two) derivatives of the $j$-function. First, we state the EC conjecture for $j$ and its derivatives. 
Let $J:\mathbb{H}\rightarrow \mathbb{C}^3$ be given by $$J: z\mapsto (j(z), j'(z), j''(z)).$$ We extend $J$ to $\mathbb{H}^n$ by defining
$$J: \bar{z}\mapsto (j(\bar{z}),j'(\bar{z}),j''(\bar{z}))$$ where $j^{(k)}(\bar{z}) = (j^{(k)}(z_1), \ldots, j^{(k)}(z_n))$ for $k=0,1,2$. Let $\Gamma_J \seq \h^n \times \C^{3n}$ be the graph of $J$ for any $n$. Note that we consider only the first two derivatives of $j$, for the higher derivatives are algebraic over those (see \cite{mahler} and also Section \ref{subsect-j-reducts}), and from a transcendence point of view adding higher derivatives would not change anything.

\begin{conjecture}[Existential Closedness for $J$]\label{J-EC}
Let $V \subseteq \h^n \times \C^{3n}$ be an irreducible $J$-broad, $J$-free and  {$\h$-free} variety defined over $\C$. Then $V\cap \Gamma_J \neq \emptyset$.
\end{conjecture}

Here $J$-\emph{broadness} and $J$-\emph{freeness} are analogues of $j$-broadness and $j$-freeness respectively (see Section \ref{subsect-broad-free} for definitions).  
Note that the aforementioned differential EC statement for the $j$-function does in fact incorporate the derivatives of $j$ and so it is a differential variant of Conjecture \ref{J-EC} (see \cite{Aslanyan-Eterovic-Kirby-Diff-EC-j} for details).

\begin{definition}
For a subgroup $G \subseteq \GL_2(\C)$ define a relation $$\B^G_J:= \left\{\left( z, j(gz), \frac{d}{dz}j(gz), \frac{d^2}{dz^2}j(gz) \right):  g \in G, gz \in \h \right\} \subseteq \C^{4}.$$
By abuse of notation for each $n$ we denote the set $$\{ (\bar{z}, \bar{w}, \bar{w}_1, \bar{w}_2) : (z_k, w_k, w_{1,k}, w_{2,k}) \in \B_J^G \mbox{ for all } k \}\seq \C^{4n}$$ by $\B_J^G$.
\end{definition}


 {Note that for trivial $G$ we have $\B_J^G = \Gamma_J$ and Conjecture \ref{J-EC} implies that if $V$ is $J$-broad, $J$-free and $\h$-free then $V\cap \B_j^G \neq \emptyset$ for any group $G$.}

We prove the following Existential Closedness statement for the blurred $J$-function.

\begin{theorem}\label{thm-ec-blur-J}
Let $V \subseteq \C^{4n}$ be an irreducible $J$-broad and $J$-free variety defined over $\C$, and let $G\subseteq \GL_2(\C)$ be a subgroup which is dense in the complex topology. Then $ V\cap \B_J^G$ is dense in $ V$ in the complex topology. In particular, if $K\subseteq \C$ is a subfield with $K\nsubseteq \R$, then $V \cap \B_J^{\GL_2(K)} \neq \emptyset$.
\end{theorem}

The proof of this result is based on the Ax-Schanuel theorem for the $j$-function \cite{Pila-Tsim-Ax-j} and the Remmert open mapping theorem from complex analytic geometry. It is an adaptation of Kirby's proof of EC for the blurred exponentiation \cite[Proposition 6.2]{Kirby-blurred}, however dealing with $j$ and its derivatives makes the proof slightly more involved.

Furthermore, we establish a stronger result for $j$ without derivatives. When $G = \GL_2^+(\Q)$ we call $\B^G_j$ the \emph{approximate} $j$-function and denote it by $\A_j$.

\begin{theorem}\label{thm-ec-appr-j}
Let $V \subseteq \h^n \times \C^{n}$ be an irreducible $j$-broad and $j$-free variety defined over $\C$. Then $V\cap \A_j$ is dense in $V$ in the complex topology.
\end{theorem}

This theorem immediately implies that $V\cap \B_j^{\GL_2(K)}$ is dense in $V$ for any subfield $K\seq \C$.

Note that the direct analogue of the \emph{approximate exponentiation} from \cite{Kirby-blurred} for the $j$-function would in fact be the relation $\B_j^{\GL_2(\Q(i))}$, and $\A_j$ is a finer blurring. EC for  $\B_j^{\GL_2(\Q(i))}$ is given by Theorem \ref{j-blur-ec-intro}, while Theorem \ref{thm-ec-appr-j} gives a stronger result. To prove it we invoke some tools from o-minimality along with basic complex geometry and the Ax-Schanuel theorem. The analogue of Theorem \ref{thm-ec-appr-j} for blurred $\exp$ is not known and our methods would fail there. This is a significant difference between this paper and \cite{Kirby-blurred}.

It is also worth noting that in fact we prove a stronger form of Theorem \ref{thm-ec-appr-j}. Namely, we show that if $G$ is dense in a certain real Lie subgroup of $\GL_2(\R)$ of dimension $2$, then the conclusion of the theorem holds for $\B_j^G$ (see Theorem \ref{thm-ec-appr-j-stronger}). Moreover, our proof will show that the theorem also holds for $G = \SL_2(\Q)$ and even for the subgroup thereof consisting of the upper triangular matrices.

Finally, at the end of the paper we show that for a suitably chosen countable algebraically closed field $C \seq \C$ the structures $\left(\C; +, \cdot, \B_j^{\GL_2(C)}\right)$ and $\left(\C; +, \cdot, \B_J^{\GL_2(C)}\right)$ are model theoretically tame. More precisely, we prove that these structures are elementarily equivalent to certain reducts of differentially closed fields, which gives a complete axiomatisation of their theories and shows that they are $\omega$-stable of Morley rank $\omega$ and near model complete (i.e. every formula is equivalent to a Boolean combination of existential formulas). We also show that they are \emph{quasiminimal}, that is, every definable set (in one variable) is either countable or co-countable.

\section{Preliminaries}\label{section-prelim}

\subsection{Modular polynomials and special varieties}\label{subsect-modular-poly}

The group $\GL_2(\C)$ acts on the Riemann sphere by linear fractional transformations and its subgroup $\GL_2^+(\R)$ -- the group of $2 \times 2$ real matrices with positive determinant -- acts on $\h$. This induces actions of subgroups of $\GL_2^+(\R)$ on $\h$. The $j$-function satisfies certain algebraic ``functional equations'' under the action of $\GL_2^+(\Q)$. More precisely, there is a countable collection of irreducible polynomials $\Phi_N\in \Z[X,Y]~ (N\geq 1)$, called \emph{modular polynomials}, such that for any $z_1, z_2\in \h$
$$\Phi_N(j(z_1),j(z_2))=0 \mbox{ for some } N \mbox{ iff } z_2=gz_1 \mbox{ for some } g\in \GL_2^+(\Q). $$

In particular, if $\tau\in \h$ is a quadratic number then $j(\tau)$ is algebraic. These numbers are known as \emph{special values} of $j$ or as \emph{singular moduli} (they are the $j$-invariants of elliptic curves with complex multiplication). 
We refer the reader to \cite{Lang-elliptic} for details.

\begin{definition}
A \emph{special} subvariety of $\mathbb{C}^n$ (with coordinates $\bar{w}$) is an irreducible component of a variety defined by modular equations, i.e. equations of the form $\Phi_N(w_k, w_l) = 0$ for some $1\leq k, l \leq n$ where $\Phi_N$ is a modular polynomial. 
\end{definition}

\subsection{Broad and free varieties}\label{subsect-broad-free}

\begin{notation} We will use the following notation.
\begin{itemize} 
    \item For a positive integer $n$ the tuple $(1,\ldots,n)$ is denoted by $(n)$, and $\bar{k}\subseteq (n)$ means that $\bar{k} = (k_1,\ldots,k_l)$ for some $1\leq k_1 < \ldots < k_l\leq n$.
    
    \item The coordinates of $\C^{2n}$ will be denoted by $(z_1,\ldots,z_n,w_1,\ldots,w_n)$. The projection on the last $n$ coordinates will be denoted by $\pr_w$. Depending on the context, we will use $\bar{z}$ or $\bar{w}$ for the coordinates of $\C^n$. 
    
    \item For a tuple $\bar{k} = (k_1,\ldots,k_l)\subseteq (n)$ define the projection map $\pi_{\bar{k}}:\C^{n} \rightarrow \C^{l}$ by
$$\pi_{\bar{k}}:(z_1,\ldots,z_n)\mapsto (z_{k_1},\ldots,z_{k_l}).$$
Further, define $\pr_{\bar{k}}:\C^{2n}\rightarrow \C^{2l}$ by
$$\pr_{\bar{k}}:(\bar{z},\bar{w})\mapsto (\pi_{\bar{k}}(\bar{z}),\pi_{\bar{k}}(\bar{w})).$$

   \item The coordinates of $\C^{4n}$ will be denoted by $(\bar{z},\bar{w},\bar{w}_1,\bar{w}_2)$, and $\ppr_{w}$ will denote the projection on the second $n$ coordinates.
    
    \item For a tuple $\bar{k} = (k_1,\ldots,k_l)\subseteq (n)$ define a map $\ppr_{\bar{k}}:\C^{4n}\rightarrow \C^{4l}$ by
$$\ppr_{\bar{k}}:(\bar{z},\bar{w},\bar{w}_1,\bar{w}_2)\mapsto (\pi_{\bar{k}}(\bar{z}),\pi_{\bar{k}}(\bar{w}),\pi_{\bar{k}}(\bar{w}_1),\pi_{\bar{k}}(\bar{w}_2)).$$
\end{itemize}
\end{notation}

\begin{definition}

\begin{itemize}

    \item[] 
    \item An algebraic variety $V \subseteq \C^{2n}$ is \emph{$j$-broad} if for any $\bar{k}\subseteq (n)$ of length $l$ we have $\dim \pr_{\bar{k}} V \geq l$. 
    
    \item An algebraic variety $V \subseteq \C^{2n}$ is \emph{$j$-free} if $\pr_w V$ {has no constant coordinate and} is not contained in any proper  special subvariety of $\C^n$.
    
    \item An algebraic variety $V \subseteq \C^{4n}$ is \emph{$J$-broad} if for any $\bar{k}\subseteq (n)$ of length $l$ we have $\dim \ppr_{\bar{k}} V \geq 3l$.
    
    \item An algebraic variety $V \subseteq \C^{4n}$ is \emph{$J$-free} if $\ppr_w V$ is not contained in any proper  special subvariety of $\C^n$, {and no coordinate is constant on the projection of $V$ to the last $3n$ coordinates}.
    
    \item  An algebraic variety $V$ in $\h^n \times \C^{n}$ or $\h^n \times \C^{3n}$ is $\h$-\emph{free} if no equation of the form $z_i = gz_k$ with $g \in \GL_2^+(\Q)$, {or of the form $z_k = c$ with $c\in \h$ a constant,} holds on $V$.
\end{itemize}

\end{definition}

\subsection{The Ax-Schanuel theorem for the $j$-function}\label{subsect-Ax-Sch}

The Ax-Schanuel theorem, proven by Pila and Tsimerman, will play a key role in the proofs of the main theorems.

\begin{fact}[Complex Ax-Schanuel for $J$, \cite{Pila-Tsim-Ax-j}]
Let $V\subseteq \mathbb{C}^{4n}$ be an algebraic variety and let $U$ be an analytic component of the intersection $V \cap \Gamma_J$. If $\dim U > \dim V - 3n$ and no coordinate is constant on $\ppr_{w} U$ then $\ppr_{w} U$ is contained in a proper special subvariety of $\mathbb{C}^n$.
\end{fact}

We will need a uniform version of this theorem. We introduce some notation first. 
For $g\in \GL_2(\mathbb{C})$ let $\mathbb{H}^g := g^{-1}\mathbb{H}$ and let $j_g: \mathbb{H}^g \rightarrow \mathbb{C}$ be the function $j_g(z) = j(gz)$. For a tuple $\bar{g} = (g_1, \ldots, g_n) \in \GL_2(\mathbb{C})^n$ let $\mathbb{H}^{\bar{g}}:= \mathbb{H}^{g_1}\times \cdots \times \mathbb{H}^{g_n}$ and define the functions
$$j_{\bar{g}}: \mathbb{H}^{\bar{g}}\rightarrow \mathbb{C}^n : (z_1,\ldots,z_n) \mapsto (j_{g_1}(z_1), \ldots, j_{g_n}(z_n))$$
and $$J_{\bar{g}} = (j_{\bar{g}},j'_{\bar{g}},j''_{\bar{g}}): \mathbb{H}^{\bar{g}}\rightarrow \mathbb{C}^{3n}: \bar{z} \mapsto (j_{\bar{g}}(\bar{z}),j'_{\bar{g}}(\bar{z}),j''_{\bar{g}}(\bar{z}))$$ where the derivation is coordinatewise and $$j'_{g_i}(z_i) =\frac{d}{dz_i} j(g_iz_i) = j'(g_iz_i)\cdot g_i'z_i,~ j_{g_i}''(z_i) = \frac{d^2}{dz_i^2}j(g_iz_i).$$ Here $j'(z) = \frac{d}{dz}j(z)$ and $g'z = \frac{d}{dz}(gz)$ for $g\in \GL_2(\mathbb{C})$. We let $\Gamma_j^{\bar{g}}\subseteq \mathbb{H}^{\bar{g}}\times \mathbb{C}^{n}$ and $\Gamma_J^{\bar{g}}\subseteq \mathbb{H}^{\bar{g}}\times \mathbb{C}^{3n}$ denote the graphs of $j_{\bar{g}}$ and $J_{\bar{g}}$ respectively.

\begin{fact}[Uniform Ax-Schanuel for $J$, {\cite[Theorem 7.7]{Aslanyan-weakMZPD}}]
Let $V_{\bar{s}}\subseteq \mathbb{C}^{4n}$ be a parametric family of algebraic varieties where $\bar{s}$ ranges over a constructible set $Q$. Then there is a finite collection $\Sigma$ of proper special subvarieties of $\mathbb{C}^n$ such that for every $\bar{s}\in Q(\C)$ and $\bar{g}\in \GL_2(\C)^n$, if $U$ is an analytic component of the intersection $V_{\bar{s}} \cap \Gamma_J^{\bar{g}}$ with $\dim U > \dim V_{\bar{s}} - 3n$, and no coordinate is constant on $\ppr_{w} U$, then $\ppr_{w} U$ is contained in some $T \in \Sigma$.
\end{fact}

This theorem implies the following statement for $j$.

\begin{fact}[Uniform Ax-Schanuel for $j$]\label{complex-Ax-Sch-uniform}
Let $(V_{\bar{s}})_{\bar{s}\in Q}$ be a parametric family of algebraic varieties in $\C^{2n}$. Then there is a finite collection $\Sigma$ of proper special subvarieties of $\C^n$ such that for every $\bar{s}\in Q(\C)$ and every $\bar{g}\in \GL_2(\mathbb{C})^n$, if $U$ is an analytic component of the intersection $V_{\bar{s}} \cap \Gamma_j^{\bar{g}}$ with $\dim U > \dim V_{\bar{s}} - n$, and no coordinate is constant on $\pr_w U$, then $\pr_w U$ is contained in some $T \in \Sigma$.
\end{fact}

\section{Proof of Theorem \ref{thm-ec-appr-j}}\label{section-proof-appr}

We will prove a slightly stronger theorem. Consider the group $$\mathcal{G}:= \begin{pmatrix}
\R^{>0} & \R \\
0 & 1
\end{pmatrix} =  \left\{ \begin{pmatrix}
a & b \\
0 & 1
\end{pmatrix} : a\in \R^{>0}, b \in \R \right\} \subseteq \GL_2^+(\R).$$
This is a two-dimensional real Lie group. 

\begin{theorem}\label{thm-ec-appr-j-stronger}
Let $V \subseteq \h^n \times \C^{n}$ be an irreducible $j$-broad and $j$-free variety defined over $\C$, and let $G\subseteq \mathcal{G}$ be a dense subgroup (in the Euclidean topology). Then $V\cap \B_j^G$ is dense in $ V$ in the complex topology.  {In particular, this holds for $G = \GL_2(\Q)\cap \mathcal{G}$.}
\end{theorem}

In the proof we will work in the o-minimal structure $\R_{\an, \exp}$ -- the expansion of the real field by restricted analytic functions and the full exponential function -- and will apply the following special case of the fibre dimension theorem in o-minimal structures. We refer the reader to \cite[Chapter 4, Corollary 1.6]{vdDries-tame-topology} for the fibre dimension theorem and to \cite{vddries-macint-marker,vddries-miller} for details on $\R_{\an,\exp}$. 

\begin{fact}\label{prop-omin-finite-fibres}
Let $(R;<,\ldots)$ be an o-minimal structure. Let also $X \subseteq R^n$ be a definable set and $f:X\rar R^m$ be a definable map with finite fibres. Then $\dim f(X) = \dim X$.
\end{fact}

In this section we will deal with real and complex dimensions of (locally) analytic sets. To distinguish between them we will use the notation $\dim_{\C}$ and $\dim_{\R}$ for the complex and real dimensions respectively.

\begin{proof}[Proof of Theorem \ref{thm-ec-appr-j-stronger}]

Let $V \subseteq \h^n \times \C^{n}$ be an irreducible $j$-broad and $j$-free variety. We may assume $\dim V=n$ by intersecting $V$ with generic hyperplanes and reducing its dimension. {See \cite[Lemma 4.31]{Aslanyan-adequate-predim} for definitions and the proof of broadness of such an intersection. Freeness is also easy to verify}.

Pick a tuple $\bar{k}=(k_1,\ldots,k_l)\subseteq (1,\ldots,n)$. For a tuple $\bar{s}\in \pr_{\bar{k}}V \subseteq \C^{2l}$ consider the fibre $V_{\bar{s}}\subseteq \C^{2(n-l)}$ above $\bar{s}$. This gives a parametric family of algebraic varieties. Let $\Sigma_{\bar{k}}$ be the collection of special subvarieties of $\C^{n-l}$ given by uniform Ax-Schanuel for this family. 

Further, by the fibre dimension theorem (see \cite[Chapter 1, \S 6, Theorem 1.25]{Shafarevich}) there is a proper Zariski closed subset $W_{\bar{k}}$ of $\pr_{\bar{k}}V$ such that if $\bar{s}\notin W_{\bar{k}}$ then 
\begin{equation}\label{eq-generic-fibre}
    \dim V_{\bar{s}} = \dim V - \dim \pr_{\bar{k}}V \leq n-l
\end{equation}
where the last inequality follows from the assumption that $V$ is $j$-broad. 

Consider the set $$V^*:=\left\{ \bar{e}\in V: \pr_{\bar{k}}\bar{e} \notin W_{\bar{k}},~ \pr_w \pr_{\bar{k}}\bar{e} \notin \bigcup_{S\in \Sigma_{\bar{k}}} S, \mbox{ for all } \bar{k} \right\}.$$
Clearly, $V^*$ is a Zariski open subset of $V$ and $V^*\neq \emptyset$ as $V$ is $j$-free. Hence $V^*$ is open and dense in $V$ with respect to the complex topology. Moreover, $V^*$ is locally analytic.


{We next observe that $\mathcal{G}$ acts transitively on $\h$. Moreover, for any two points $z_1, z_2 \in \h$ there is a unique element of $\mathcal{G}$ that maps $z_1$ to $z_2$. To this end, let $z_1=x+iy$ and $z_2 = u+iv$ where $x, u \in \R,~  y, v \in \R^{>0}$. Then it is straightforward to check that the matrix $$g(z_1,z_2):= \begin{pmatrix}
\frac{v}{y} & u- \frac{xv}{y} \\
0 & 1
\end{pmatrix} \in \mathcal{G}$$ maps $z_1$ to $z_2$, and simple calculations show that it is the only element of $\mathcal{G}$ with that property. }

 {Pick a fundamental domain $\F \seq \h$ of the action of $\SL_2(\Z)$. Since the restriction $j|_{\F}:\F \rar \C$ is bijective, there is an inverse function $j^{-1}:\C \rar \F$.} Now define a map 
\begin{gather*}
    \theta: \h^n \times \C^n \rightarrow \mathcal{G}^n,\\ 
    \theta: (\bar{z},\bar{w})\mapsto (g(z_1,j^{-1}(w_1)),\ldots,g(z_n,j^{-1}(w_n))),
\end{gather*}
{and let $\theta^*:=\theta|_{V^*}$ be the restriction of $\theta$ to $V^*$.}

It is well known that the restriction of $j$ to any fundamental domain is definable in the o-minimal structure $\R_{\an,\exp}$ (see \cite{Peterzil-Starchenko-unifrom-def}). Hence, $\theta$ is a definable map. Since $\mathcal{G}$ and $V^*$ are definable in the field structure of $\R$, the map $\theta^*$ is definable in $\R_{\an,\exp}$.

\begin{claim}
All fibres of $\theta^*$  are finite.
\end{claim}
\begin{proof}
{The main idea of the proof is that even though $\theta^*$ is not a complex analytic map (it is of course real analytic), its fibres are contained in complex analytic sets of dimension $0$.} Pick a tuple $\bar{g}\in \mathcal{G}^n$. It is easy to see that $$\theta^{-1}(\bar{g})\seq  \Gamma_j^{\bar{g}}.$$ Indeed,
if $\theta(\bar{z},\bar{w}) = \bar{g}$ then $g_kz_k = j^{-1}(w_k) \mbox{ for all } k$, hence $w_k = j(g_kz_k) \mbox{ for all } k.$ Thus, $$(\theta^*)^{-1} (\bar{g}) \seq V^* \cap \Gamma_j^{\bar{g}}.$$ 
Observe that $V^*\cap \Gamma_j^{\bar{g}}$ is a locally analytic subset of $\h^n \times \C^n$. Let $U$ be a non-empty irreducible component of it. We will show that $\dim_{\C} U =0$. Pick a point $\bar{e}\in U$. Note that a coordinate $z_k$ is constant on $U$ if and only if the corresponding coordinate $w_k$ is also constant on $U$, i.e. constant coordinates on $U$ come in pairs. Let $z_{k_1},\ldots,z_{k_l}$ and the corresponding $w$-coordinates list all of the constant coordinates on $U$. Set $\bar{k}:=(k_1,\ldots,k_l)$ and $\bar{t}:=\pr_{\bar{k}}\bar{e}$, and consider the fibre $U_{\bar{t}}$. Then $\dim_{\C} U_{\bar{t}} = \dim_{\C} U$ and $U_{\bar{t}}$ has no constant coordinates. Moreover, $U_{\bar{t}}$ is an analytic component of $V^*_{\bar{t}}\cap \Gamma_j^{\bar{h}}$ {where $\bar{h}$ is the tuple of all $g_p$ with $p\neq k_1,\ldots,k_l$}. Since $\bar{e}\notin S$ for any $S\in \Sigma_{\bar{k}}$, we conclude by the Ax-Schanuel theorem that $\dim_{\C} U_{\bar{t}} = \dim_{\C} V^*_{\bar{t}} - (n-l) \leq 0$ by \eqref{eq-generic-fibre}, hence $\dim_{\C} U_{\bar{t}} = \dim_{\C} U =0$.

Now since $U$ is locally analytic, connected and of dimension $0$, it must be a singleton. We showed that all analytic components of the set $V^* \cap \Gamma_{\bar{g}}$ are singletons, hence $V^* \cap \Gamma_j^{\bar{g}}$ is discrete. It must actually be finite, for it is definable in $\R_{\an, \exp}$.
\end{proof}

Thus, the map $\theta^* : V^* \rightarrow \mathcal{G}^n$ has finite fibres and is definable in the o-minimal structure $\R_{\an, \exp}$. The real dimension of $V^*$ is clearly $2n$, as is the dimension of $\mathcal{G}^n$. Since the fibres of $\theta^*$ are finite, the image $\theta^*(V^*)$ must have dimension $2n$ by Fact \ref{prop-omin-finite-fibres}. In particular, $\theta^*(V^*)$ must contain a cell of dimension $2n$ which means that it has non-empty interior. Therefore, $\theta^*(V^*)$ has non-empty intersection with any dense subset of $\mathcal{G}^n$, in particular with $G^n$. 

Moreover, if $O\subseteq V^*(\C)$ is an open subset (definable in $\R$), then $\dim_{\R}O = \dim_{\R}V^* = 2n$ and $\dim_{\R} \theta^*(O) = 2n$. So $\theta^*(O)\subseteq \mathcal{G}^n$ has non-empty interior and $\theta^*(O)\cap G^n \neq \emptyset$. Therefore, $O \cap \B_j^G \neq \emptyset$. This means that $V^*\cap \B_j^G$ is dense in $V^*$. Since $V^*$ is dense in $V$, we conclude that $V\cap \B_j^G$ is dense in $V$.
\end{proof}

\begin{remark}
Slightly adapting the proof shows that in Theorem \ref{thm-ec-appr-j-stronger} we can replace $\mathcal{G}$ by any Lie subgroup of $\GL^+_2(\R)$ that acts transitively on $\h$, such as the subgroup of $\SL_2(\R)$ consisting of upper triangular matrices. 
\end{remark}

\section{Proof of Theorem \ref{thm-ec-blur-J}}\label{section-proof-blur-J}
\setcounter{equation}{0}

First, we observe that the relation $\B_J^G$ can be expressed in terms of the functions $J_{\bar{g}}$ and their graphs $\Gamma_J^{\bar{g}}$ (defined in Section \ref{subsect-Ax-Sch}) as follows. For a subgroup $G \subseteq \GL_2(\C)$ we have $$\B^G_J= \{(\bar{z}, J_{\bar{g}}(\bar{z})): \bar{g}\in G^n, \bar{z}\in \h^{\bar{g}} \} = \bigcup_{\bar{g}\in G^n} \Gamma_J^{\bar{g}} \subseteq \C^{4n}.$$

Now let $V \subseteq \C^{4n}$ be $J$-broad and $J$-free. As in the proof of Theorem \ref{thm-ec-appr-j-stronger}, we may assume $\dim V =3n$. We shall define a Zariski open subset $V'$ of $V$ as in the previous section. For a tuple $\bar{k}=(k_1,\ldots,k_l)\subseteq (1,\ldots,n)$ consider the parametric family of the fibres $V_{\bar{s}}\subseteq \C^{4(n-l)}$ for $\bar{s}\in \ppr_{\bar{k}}V \subseteq \C^{4l}$. Let $\Sigma_{\bar{k}}$ be the collection of special subvarieties of $\C^{n-l}$ given by uniform Ax-Schanuel with derivatives for this family. 
Let also $W_{\bar{k}}\seq \ppr_{\bar{k}}V$ be a Zariski closed subset such that 
\begin{equation*}\label{eq-generic-fibre-deriv}
    \dim V_{\bar{s}} = \dim V - \dim \ppr_{\bar{k}}V \leq 3(n-l) \mbox{ whenever } \bar{s}\notin W_{\bar{k}}.
\end{equation*}
Finally, let $$V':=V^{\reg} \cap \left\{ \bar{e}\in V: \ppr_{\bar{k}}\bar{e} \notin W_{\bar{k}},~ \ppr_w \ppr_{\bar{k}}\bar{e} \notin \bigcup_{S\in \Sigma_{\bar{k}}} S, \mbox{ for all } \bar{k} \right\},$$
{where $V^{\reg}$ denotes the set of regular points of $V$. Then $V'$ is a smooth Zariski open subset of $V$.}



\begin{lemma}
Given a complex number $D\neq 0$, there are an open and dense subset $U \subseteq  \C^4$ (in the complex topology) and an analytic map $\zeta_D: U \rightarrow \GL_2(\C)$ such that for any $(z, w, w_1, w_2)\in U$ if $g=\zeta_D(z,w,w_1,w_2)$ then we have $\det g=D$ and 
\begin{equation}\label{eq-J}
    w=j(gz), w_1=(j(gz))', w_2=(j(gz))''.
\end{equation}
\end{lemma}
\begin{proof}
We solve the system of equations \eqref{eq-J} with respect to $g =\begin{pmatrix}
a & b \\
c & d
\end{pmatrix}.$ Observe that if $\det g = D$ then
$$(j(gz))' = j'(gz)\cdot (gz)' = j'(gz)\cdot \frac{D}{(cz+d)^2}$$ and 
$$(j(gz))''  = j''(gz)\cdot \frac{D^2}{(cz+d)^4} - j'(gz)\cdot \frac{2cD}{(cz+d)^3}.$$

Let $\F\seq \h$ be a fundamental domain for the action of $\SL_2(\Z)$. Then the map $j^{-1}:\C\rightarrow \F$ is well defined. Let $(z,w,w_1,w_2)\in \C^4$ {be such that $w_1\neq 0$ and $j'(\tilde{z}) \neq 0$ where} $\tilde{z}:=j^{-1}(w)\in \F$. Consider the equations
$$w_1 = j'(\tilde{z})\cdot \frac{D}{(cz+d)^2},~ w_2 = j''(\tilde{z})\cdot \frac{D^2}{(cz+d)^4} - j'(\tilde{z})\cdot \frac{2cD}{(cz+d)^3}.$$ 

Solving these with respect to $c$ and $d$ we get
\begin{gather*}
    c = \frac{1}{2w_1}\cdot \sqrt{\frac{Dj'(\tilde{z})}{w_1}} \cdot \left[ j''(\tilde{z})\cdot \left( \frac{w_1}{j'(\tilde{z})}\right)^2 - w_2 \right],\\
    d = \sqrt{\frac{Dj'(\tilde{z})}{w_1}} - cz.
\end{gather*}

Further, we find $a$ and $b$ from the equations $$\frac{az+b}{cz+d} = \tilde{z},~ ad-bc=D.$$

The solutions for $a,b,c,d$ are functions of $z,w,w_1,w_2$, and it is clear that \eqref{eq-J} holds for those functions. Let $\F^0$ be the interior of $\F$ and let $\C':=j(\F^0)$. Then $\C'$ is open and dense in $\C$ and $j^{-1}:\C' \rightarrow \F^0$ is holomorphic. Choosing a branch of the square root we see that $a,b,c,d$ are holomorphic on an open dense subset $U \subseteq \C \times \C' \times \C^2$. {Note that we also choose $U$ such that $j'(j^{-1}(w))\neq 0$ and $w_1\neq 0$ on $U$.} 
\end{proof}

Now define a map $$\theta: = (\zeta_1, \ldots, \zeta_1):U^n \rightarrow \SL_2(\C)^n.$$ Consider the set $V^*:=V'\cap U^n$ and the restriction $\theta^*:=\theta|_{V^*}$. {The coordinates $w$ and $w_1$ are non-constant on $V$ by $J$-freeness and $V'$ is Zariski dense in $V$, hence $V^*$ is of dimension $3n$}. It is also evident that $V^*$ is a smooth analytic subset of $U^n$ and $\theta^*$ is a holomorphic map. Given a tuple of matrices $\bar{g} \in \theta^*(V^*)$, we have $$(\theta^*)^{-1}(\bar{g}) \subseteq V'\cap \Gamma_J^{\bar{g}}.$$

As in the proof of the claim in the previous section we can employ the Ax-Schanuel theorem to show that $V'\cap \Gamma_J^{\bar{g}}$ is discrete. Thus, $\theta^*$ has discrete fibres. Since $\dim_{\C} V^* = \dim_{\C} \SL_2(\C)^n = 3n$, by Remmert's open mapping theorem (see \cite[\S V.6, Theorem 2]{Loj-analytic-geom}) we conclude that $\theta^*$ is an open map. Therefore for any open subset $O\subseteq V^*$ the image $\theta^*(O)\subseteq \SL_2(\C)^n$ is open, so it intersects $G^n$. Thus, any open subset of $V^*$ intersects $\B_J^G$ and so $V\cap \B_J^G$ is dense in $V$.

Finally, if $C$ is a subfield of $\C$ not contained in $\R$, then $C$ is dense in $\C$ and $\SL_2(C)$ is dense in $\SL_2(\C)$ and the second part of Theorem \ref{thm-ec-blur-J} follows.\\

We can deduce Theorem \ref{j-blur-ec-intro} from Theorem \ref{thm-ec-blur-J}.

\begin{proof}[Proof of Theorem \ref{j-blur-ec-intro}]
If $V \subseteq \C^{2n}$ is $j$-broad and $j$-free then $\tilde{V}:= V \times \C^{2n}\seq \C^{4n}$ is $J$-broad and $J$-free, hence $\tilde{V}(\C) \cap \B_J^G$ is dense in $\tilde{V}(\C)$. Therefore, $ V\cap \B_j^G$ is dense in $ V$.
\end{proof}

\section{Model theoretic properties of the blurred $j$-function}\label{section-model-theoretic-properties}

In this section we show that for a suitably chosen countable algebraically closed subfield $C\seq \C$ and the group $G:=\GL_2(C)$ the structures $\C_{\B_j^G}:=(\C; +, \cdot, \B_j^{G})$ and $\C_{\B_J^G}:=(\C; +, \cdot, \B_J^{G})$ are model theoretically tame, where $\B_j^{G}$ and $\B_J^G$ are considered respectively as binary and $4$-ary relations. In particular, we get a first-order axiomatisation of those structures and show that they are $\omega$-stable of Morley rank $\omega$.

\subsection{$j$-reducts of differentially closed fields}\label{subsect-j-reducts}

The $j$-function satisfies a third order algebraic differential equation over $\mathbb{Q}$ \cite{mahler}. Namely {(see e.g. \cite[p. 20]{Masser-heights})}, $\Psi(j,j',j'',j''')=0$ where 
$$\Psi(w,w_1,w_2,w_3)=\frac{w_3}{w_1}-\frac{3}{2}\left( \frac{w_2}{w_1} \right)^2 + \frac{w^2-1968w+2654208}{2w^2(w-1728)^2}\cdot w_1^2.$$

Thus $$\Psi(w,w',w'',w''')=Sw+R(w)(w')^2,$$
where $S$ denotes the \emph{Schwarzian derivative} defined by $$Sw = \frac{w'''}{w'} - \frac{3}{2} \left( \frac{w''}{w'} \right) ^2$$ and $$R(w)=\frac{w^2-1968w+2654208}{2w^2(w-1728)^2}.$$ Note that all functions $j(gz)$ with $g \in \GL_2(\mathbb{C})$ satisfy the differential equation $\Psi(w,w',w'',w''')=0$ and in fact all solutions are of that form (see \cite[Lemma 4.2]{Freitag-Scanlon} or \cite[Lemma 4.1]{Aslanyan-adequate-predim}).

In a differential field $(K;+, \cdot, ')$ for a non-constant $z \in K$ define a derivation $\partial_z:K \rightarrow K$ by $\partial_z: w\mapsto \frac{w'}{z'}$. The expression $\Psi(w,\partial_zw,\partial^2_zw,\partial^3_zw)$ is a differential rational function in $z$ and $w$, and we let $\chi(z,w)$ be its numerator. Let also $E_j\subseteq K^2$ be a binary relation that is interpreted as the set of solutions of the equation $\chi(z,w)=0$. The reduct $(K;+,\cdot, E_j)$ will be denoted by $K_{E_j}$.

\begin{fact}[{\cite[Theorem 4.39]{Aslanyan-adequate-predim}} and {\cite[Theorem 1.1]{Aslanyan-Eterovic-Kirby-Diff-EC-j}}]
Let $(K;+,\cdot, D)$ be a differentially closed field with field of constants $C$. Then $T_j:=\Th(K_{E_j})$ is axiomatised by the following axioms and axiom schemes.

\begin{enumerate}
\item[\rm A1] $K$ is an algebraically closed field of characteristic $0$. 

\item[\rm A2] $C:=\{ c \in K: E_j(0,c)\}$ is an algebraically closed subfield. Further, $C^2 \subseteq E_j(K)$ and if $\left(z,j\right) \in E_j(K)$ and one of $z,j$ is constant then both of them are constants.

\item[\rm A3] If $\left(z,w\right) \in E_j$ then for any  $g \in \SL_2(C)$, 
$\left(gz,w\right) \in E_j.$
Conversely, if for some $w$ we have $\left(z_1,w\right), \left(z_2,w\right) \in E_j$ then $z_2=gz_1$ for some $g \in \SL_2(C)$.

\item[\rm A4] If $\left(z,w_1\right) \in E_j$ and $\Phi_N(w_1,w_2)=0$ for some $w_2$ and some modular polynomial $\Phi_N(X,Y)$ then $\left(z,w_2\right) \in E_j$.

\item[\rm AS] If $\left(z_i,w_i\right) \in E_j,~ i=1,\ldots,n,$ with
$$\td_CC\left(\bar{z},\bar{w}\right) \leq n,$$ then $\Phi_N(w_i,w_k)=0$ for some $N$ and some $1\leq i < k \leq n$, or $w_i \in C$ for some $i$.

\item[\rm EC] For each $j$-broad variety $V \subseteq K^{2n}$  the intersection $E_j(K) \cap V(K)$ is non-empty.

\item[\rm NT] There is a non-constant element in $K$.
\end{enumerate}

Furthermore, $T_j$ is $\omega$-stable of Morley rank $\omega$, and near model complete, i.e. every formula is equivalent to a Boolean combination of existential formulas modulo $T_j$.
\end{fact}

Here AS should be understood as the uniform Ax-Schanuel which is first-order expressible (see \cite[p. 26]{Aslanyan-adequate-predim} and \cite[Theorem 4.6]{Aslanyan-weakMZPD}).

Similarly, if we let $E_J$ be the set of all $4$-tuples $(w,\partial_zw,\partial^2_zw,\partial^3_zw)\in K^4$ with $\chi(z,w)=0$ then the first-order theory of $K_{E_J}:=(K;+,\cdot, E_J)$ is well understood. In particular, a complete axiomatisation has been obtained in \cite[Theorem 5.20]{Aslanyan-adequate-predim}.

\subsection{First-order theory of the blurred $j$-function}\label{subsect-first-order-th-blur-j}

\begin{definition}[{\cite[\S 5]{Eterovic-Schan-for-j}}]
A $j$-derivation on the field of complex numbers is a derivation $\delta:\C\rar \C$ such that for any $z\in \h$ we have
$$\delta j(z) = j'(z)\delta(z),~ \delta j'(z) = j''(z)\delta(z),~ \delta j''(z) = j'''(z)\delta(z).$$
The space of $j$-derivations is denoted by $j\rm{Der}(\C)$.
\end{definition}

Let $$C:=\bigcap_{\delta \in j\rm{Der}(\C)} \ker \delta.$$ Then $C$ is a countable algebraically closed subfield of $\C$ and $j(C\cap \h)=C$ (see \cite[\S 5]{Eterovic-Schan-for-j}).

\begin{theorem}[cf. {\cite[Theorem 1.4]{Kirby-blurred}}]\label{thm-blurred-diff-reduct}
Let $C$ be as above and $G=\GL_2(C)$. Then $\Th(\C_{\B_j^G}) = T_j$. In particular, $\Th(\C_{\B_j^G})$ is $\omega$-stable with Morley rank $\omega$ and is near model complete.
\end{theorem}
\begin{proof}
Since $T_j$ is complete, it suffices to prove that $\C_{\B_j^G}\models T_j$. Axioms A1-A4 and NT are straightforward to check, AS follows from \cite[Proposition 6.2]{Eterovic-Schan-for-j} and EC follows from Theorem \ref{thm-ec-appr-j}, for it is proven in \cite[\S 4.8]{Aslanyan-adequate-predim} that in EC one may assume the variety is free.
\end{proof}

Similarly, we get the following result for the blurred $J$-function.

\begin{theorem}
Let $C$ be as above and $G=\GL_2(C)$. Then $\C_{\B_J^G}$ is elementarily equivalent to $E_J$-reducts of differentially closed fields, axiomatised by \rm{A1$'$-A4$'$,AS$'$,EC$'$,NT$'$} from \cite[\S 5]{Aslanyan-adequate-predim}. In particular, $\Th(\C_{\B_J^G})$ is $\omega$-stable with Morley rank $\omega$ and is near model complete.
\end{theorem}
\begin{proof}
This can be proven exactly as Theorem \ref{thm-blurred-diff-reduct}.
\end{proof}

\subsection{Remarks on quasiminimality}\label{subsect-quasimin}
It can be shown that the structures $\C_{\B_j^G}$ and $\C_{\B_J^G}$ are \emph{quasiminimal}, that is, every first-order definable set in either of those structures is countable or its complement is countable. We briefly explain how this can be proven in the same way as for blurred complex exponentiation. We focus on $\C_{\B_J^G}$, the case of  $\C_{\B_j^G}$ being analogous. There are three main ingredients in the proof.
\begin{itemize}
    \item[(a)] There is a natural pregeometry on $\C_{\B_J^G}$ given by the Ax-Schanuel inequality, denoted by $\cl$. See \cite[Definition 2.15]{Aslanyan-adequate-predim}.
    
    \item[(b)] $\C_{\B_J^G}$ has the \emph{countable closure property}, i.e. the $\cl$-closure of a finite set is countable.
    
    \item[(c)] $\C_{\B_J^G}$ realises a unique generic type over countable closed subsets: if $A$ is a countable $\cl$-closed subset and $a_1, a_2 \in \C \setminus A$ then $\tp(a_1/A) = \tp(a_2/A)$.
\end{itemize}

For details on (a) we refer the reader to \cite{Aslanyan-adequate-predim}, (b) follows from the proof of Theorem \ref{thm-ec-blur-J}, and (c) can be proven by constructing a back-and-forth system of partial isomorphisms from $a_1$ to $a_2$ over $A$ as in the proof of \cite[Proposition 4.38]{Aslanyan-adequate-predim} (and it relies on a ``generic'' version of Theorem \ref{thm-ec-blur-J}). Now, given a  set $X\seq \C$, definable in $\C_{\B_J^G}$, let $A$ be the closure of the finitely many parameters used in the definition of $X$. Then $A$ is countable and by (c) either $X\seq A$ or $X\supseteq \C \setminus A$.

In \cite{Kirby-blurred} Kirby proves a stronger quasiminimality result for the blurred complex exponentiation, namely, every subset which is invariant under all automorphisms is either countable or co-countable. It is likely that the methods of \cite{Kirby-blurred} will go through for the $j$-function and we will have similar stronger results for  $\C_{\B_j^G}$ and $\C_{\B_J^G}$. 

Unlike complex exponentiation, where quasiminimality is an open question, the complex field equipped with the graph of the $j$-function cannot be quasiminimal, for $\h$ is definable. However, it would be interesting to understand which blurrings of the $j$-function are quasiminimal. Since the action of $\GL_2(\C)$ factors through $\PGL_2(\C)$, we ask the following question.

\begin{question}
For which proper subgroups $G$ of $\PGL_2(\C)$ are the structures $\C_{\B_j^G}$ and $\C_{\B_J^G}$ quasiminimal?
\end{question}

One can see immediately that there are some trivial examples where $\C_{\B_j^G}$ is not quasiminimal. First, when $G$ is uncountable, the fibres of $\B_j^G$ above the second coordinate are uncountable with an uncountable complement, hence $G$ must be at most countable. Second, if $G\seq \PGL_2(\R)$, then the projection of $\B_j^G$ on the first coordinate is $\h$, therefore $\C_{\B_j^G}$ is not quasiminimal. Further, when $G$ is finite then the fibres of $\B_j^G$ above the first coordinate may be finite and of different cardinalities which allows one to define an uncountable set whose complement is also uncountable. For example, when $G$ is the group generated by $\begin{pmatrix} i & 0 \\ 0 & 1 \end{pmatrix}$, then the set $S:=\{ z\in \C: \exists ! w \B_j^G(z,w) \}$ contains $\R \cup i\R \setminus \{ 0 \}$, so it is uncountable and one can easily see that $S \setminus (\R \cup i\R)$ is at most countable, hence $\C \setminus S$ is uncountable.


It seems plausible that the above question has an affirmative answer if and only if $G \nsubseteq \PGL_2(\R)$ and $G$ is countably infinite.

\section*{Acknowledgements}
We thank the referee for useful comments that helped us improve the presentation of the paper.

\bibliographystyle {alpha}
\bibliography {ref}

\end{document}